\documentclass[12pt,letterpaper]{ISSFD_v01}
\usepackage[left=1in,right=1in,top=1in,bottom=1in]{geometry}
\usepackage{graphicx}
\usepackage{amsmath}

\usepackage{times} % Use times font
\usepackage{mathptmx} % Use same font for math mode
\DeclareSymbolFont{largesymbols}{OMX}{cmex}{m}{n} % Undo the definition of large symbols from the mathptmx package
\usepackage{authblk} % Use authblk package for author formatting
\usepackage[pdffitwindow=false,pdfstartview={FitH},colorlinks=true,citecolor=black,linkcolor=black,urlcolor=black]{hyperref}
\usepackage{microtype} % Makes detailed spacing look nicer. This is only optional -- not required.

\usepackage{xspace}

\title{Review of Lambert's problem}

\author[(1)]{David de la Torre Sangr\`a}
\author[(2)]{Elena Fantino}
\affil[(1)]{Polytechnic University of Catalonia (UPC), E.T.S.E.I.A.T. calle Colom 11, 08222 Terrassa (Spain), david.de.la.torre.sangra@upc.edu}
\affil[(2)]{Space Studies Institute of Catalonia (IEEC), Polytechnic University of Catalonia (UPC), E.T.S.E.I.A.T., Colom 11, 08222 Terrassa (Spain), elena.fantino@upc.edu}
\date{} % Don't display a date

\begin{document}
\maketitle

\begin{abstract}
Lambert's problem is the orbital boundary-value problem constrained by two points and elapsed time. It is one of the most extensively studied problems in celestial mechanics and astrodynamics, and, as such, it has always attracted the interest of mathematicians and engineers. Its solution lies at the base of algorithms for, e.g., orbit determination, orbit design (mission planning), space rendezvous and interception, space debris correlation, missile and spacecraft targeting.
There is abundancy of literature discussing various approaches developed over the years to solve Lambert's problem. We have collected more than 70 papers and, of course, the issue is treated in most astrodynamics and celestial mechanics textbooks. From our analysis of the documents, we have been able to identify six or seven main solution methods, each associated to a number of revisions and variations, and many, so to say, secondary research lines with little or no posterior development. We have ascertained plenty of literature with proposed solutions, in many cases supplemented by performance comparisons with other methods.
We have reviewed and organized the existing bibliography on Lambert's problem and we have 
performed a quantitative comparison among the existing methods for its solution. The analysis is based on the following issues: choice of the free parameter, number of iterations,generality of the mathematical formulation, limits of applicability (degeneracies, domain of the parameter, special cases and peculiarities), accuracy, and suitability to automatic execution. Eventually we have tested the performance of each code.
The solvers that incorporate the best qualities are Bate's algorithm via universal variables with Newton-Raphson and Izzo's Householder algorithm. The former is the fastest, the latter exhibits the best ratio between speed, robustness and accuracy.
\end{abstract}


\section{Introduction}
\begin{figure}[ht]
\centering
\includegraphics[width=5.0cm]{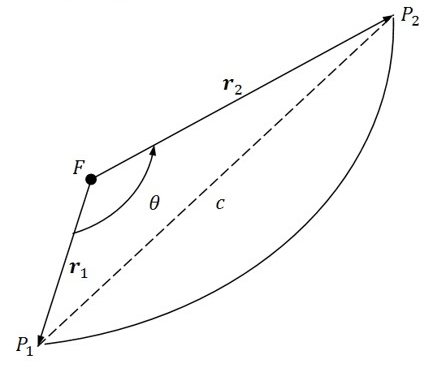}
\caption{Definition of Lambert's problem.}
\label{fig:lambert_def}
\end{figure}

This contribution deals with a prominent piece of Celestial Mechanics and Astrodynamics, the Two-Body Orbital Boundary-Value Problem (TBOBVP), also known as the Gauss' or Lambert's Problem. It is the problem of determining the Keplerian orbit connecting two positions in a given time. Figure~\ref{fig:lambert_def} illustrates its definition: $F$ is the primary (the Sun in the case of interplanetary orbits, the Earth in the context of geocentric motion), origin of the reference frame, while ${\bf r}_1$ and ${\bf r}_2$ are the positions occupied by the secondary at times $t_1$ and $t_2$, respectively (with $t_2 > t_1$), thus being $\Delta t = t_2 - t_1$ the flight time. The angle $\theta$ between ${\bf r}_1$ and ${\bf r}_2$ indicates the direction of motion. Finding the conic section that connects the two positions in the given time is equivalent to determining the velocity ${\bf v}_1$ to be imparted to the secondary at ${\bf r}_1$ to execute the transfer. The case shown as an example in Fig.~\ref{fig:lambert_def}  is an arc of an ellipse with one focus at $F$.
The TBOBVP problem arose towards the end of the 18$^{th}$ century in connection with the determination of the orbits of celestial bodies from observation.
Carl Friedrich Gauss (1777-1855) used three observation times instead of two (and with this method he was able to correctly determine the orbit of the newly-discovered Ceres). Gauss understood the potential of the TBOBVP and his developments on this subject were remarkable, with the result that the problem now bears his name. The work of Gauss merges with the previous findings of Johann Heinrich Lambert (1728-1777) who deduced the homonymous theorem:  
{\it Given the gravitational parameter $\mu= GM$, the time $\Delta t$ required to accomplish a given transfer is a function of the semimajor axis $a$ of the orbit, the sum $r_1+r_2$ of the distances from the primary at the beginning and at the end of the transfer and the length $c$ of the chord that connects such positions}, i.e.:
\begin{equation}
\sqrt{\mu} \Delta t = f(a,r_1+r_2,c). 
\label{eq:transfer_time}
\end{equation}
When the geometry of the radius vectors is fixed, there is only one free parameter left that wholly defines the transfer time between ${\bf r}_1$ and ${\bf r}_2$. In the original formulation, such parameter is the unknown semimajor axis $a$. In general, once a suitable parametrization of the orbit passing trough the given points is obtained, it is possible to rephrase Lambert's problem in terms of the evaluation of the parameter's value, such that the corresponding orbit is characterized by a transfer time that matches exactly the prescribed one. 
There exist several ways to express $f$, but no analytical closed-form solution of the transfer time equation (Eq.~\ref{eq:transfer_time}) is possible. Especially valuable and useful are the unified forms, i.e., providing one formulation valid for all three types of conic sections. 

In modern Astrodynamics, Lambert's problem has direct application in the solution of intercept and rendezvous, ballistic missiles targeting and interplanetary trajectory design. The relevance of this issue is confirmed by the vast related literature. Several authors since the time of Gauss have studied alternative formulations of the transfer time equation, adopting a variety of independent variables (free parameter) and solution methods. Due to the transcendental character of the transfer time equation, all the available approaches are based on a numerical procedure in which the value of the free parameter is searched iteratively. Then, the orbital elements or the state vector at departure and arrival are obtained by means of orbital mechanics relations. 

In this work, we present a review of the algorithms for the solution of Lambert's problem, based on quantitative comparisons. We start by discussing the relationships among the many methods, and we indicate the main representatives of the several research lines, as we identified them from our study of the literature (Sect.~\ref{sec:methods}). Then, in Sect.~\ref{sec:performance} we illustrate the tests that we carried out in order to compare the selected methods in terms of performance. We discuss the results and we draw conclusions in Sect.~\ref{sec:conclusions}.

\section{Solution methods: comparison and discussion}
\label{sec:methods}
The first formulations of Lambert's problem are due to Lagrange \cite{Lagrange1788} and Gauss \cite{Gauss1857}. However, studies on the subject turned intensive in the mid 1960s. As of today, more than 60 authors have proposed formulations and solutions to the problem. All these methods can be grouped into a number of major lines of research on the basis of the free parameter adopted. Here we highlight the most productive ones, and for each we indicate the progenitor algorithm and the successive improvements upon it:
\begin{itemize}
\item Universal variables:
\begin{itemize}
\item Lancaster \& Blanchard \cite{LANCASTER1969}, Gooding \cite{Gooding1988}~\cite{Gooding1990}, Izzo \cite{Izzo2015}.
\item Bate \cite{Bate1971}, Vallado \cite{Vallado1997}, Luo \cite{Luo2011}, Thomson \cite{Thompson2011}, Arora \cite{Arora2013a}.
\item Battin \cite{BATTIN1983}~\cite{BATTIN1984}, Loechler \cite{Loechler1988}, Shen \cite{Shen2004}, MacLellan \cite{Maclellan2005}.
\end{itemize}
\item Semi-major axis:
\begin{itemize}
\item Lagrange \cite{Lagrange1788}, Thorne \cite{Thorne1995}~\cite{Thorne2014}, Prussing \cite{Prussing2000}, Chen \cite{Chen2013}, Wailliez  \cite{Wailliez2014}.
\end{itemize}
\item Semi-latus rectum ($p$-iteration):
\begin{itemize}
\item Herrick-Liu \cite{Herrick1959}, Boltz \cite{Boltz1984}.
\item Bate \cite{Bate1971}.
\end{itemize}
\item Eccentricity vector:
\begin{itemize}
\item Avanzini \cite{Avanzini2008}, He \cite{He2010}, Zhang \cite{Zhang2010}~\cite{Zhang2011}, Wen \cite{Wen2014a}.
\end{itemize}
\item Kustaanheimo-Stiefel (K-S) regularized coordinates:
\begin{itemize}
\item Sim\'o \cite{Simo1973}
\item Kriz \cite{Kriz1976}, Jezewsky \cite{Jezewski1976}.
\end{itemize}
\end{itemize}
The work of Sun \cite{Sun1979} is also worth mentioning, being the first to extensively investigate and analyse the multi-revolution Lambert's problem and derive an expression for the minimum transfer time in this case.

We have selected representatives of each line of work and we have studied them. \cite{Lagrange1788} and \cite{Gauss1857} have been retained for historical reasons. From the Universal variables line, the selection includes Bate \cite{Bate1971}, Battin \cite{BATTIN1984}), and the work by Gooding \cite{Gooding1990} and Izzo \cite{Izzo2015}, both superseding the original algorithm of Lancaster \& Blanchard. Finally, the work of Sim\'o has been selected as representative of the K-S branch. 
In the remainder of this section, we will describe these algorithms. For reasons of space, the reader is referred to the original publications for details and more quantitative aspects.
% decir en las conclusiones que los restantes metodos se analizaran en otra publicacion

\subsection{Generic Lambert procedure}
In general, the procedure to be adopted to solve Lambert's problem consists in:
\begin{enumerate}
  \item computing the geometric parameters of the transfer;
  \item obtaining an initial guess for the free parameter;
  \item iterating on the transfer time equation until convergence;
  \item computing the velocity vectors.
\end{enumerate}

\subsection{Transfer time equation}
\label{sec:tof}
%Each Lambert algorithm relies on a particular form of the time equation (Eq.~\ref{eq:transfer_time}). The time equation from each of the selected algorithms are presented below, along with an algebraic analysis on their free parameters, singularities, and particularities.
The algorithm by Lagrange expresses the transfer time as a function of the semi-major axis $a$. This choice is inconvenient for two reasons: the solution contains a pair of conjugate orbits (i.e., the solution is not unique), and the derivative of the transfer time is singular when $a$ corresponds to the minimum-energy orbit. To overcome this difficulty, Battin \cite{Battin1999} expresses the time equation as a function of the universal variable $x$ (DEFINE) which is valid for all types of conic sections, thus obtaining a single-valued and monotonic function for the transfer time.

In his Theoria Motus \cite{Gauss1857}, Gauss develops a method to solve Lambert's problem with a system of two equations and using the sector area-to-triangle area ratio $y$. Then, Bate \cite{Bate1971} expands the original Gauss definition to any type of conic section using a power series expansion developed by Moulton \cite{Moulton1984}.

In his book, Bate \cite{Bate1971} starts from the Gauss $f$ and $g$ functions and expresses the transfer time in terms of the universal variable $z$ using Stumpff functions. $z$ is defined as $\Delta E^2$ for elliptical orbits and as $-\Delta F^2$ for hyperbolic orbits, being $\Delta E$ and $\Delta F$ the eccentric anomaly in the ellipse and in the hyperbola, respectively. Note that Bate's time equation is singular at $z=(2n\pi)^2$ with $n=1,2,3...$, and $[2n\pi]^2$ and $[2(n+1)\pi]^2$ are the limits of the domain of $n$-revolution elliptical orbits. Note also that there is a lower limit for $z$ when the transfer angle is less than $\pi$. Such limit corresponds to a hyperbolic orbit with transfer time $\Delta t = 0$. Lower values of $z$ imply negative values for the transfer time.

Sim\'o \cite{Simo1973} derives the transfer time in regularized space by means of the Levi-Civita coordinate transformation. His free parameter $z$ has a universal definition, valid for the three types of conic section. The transfer time is singular at $z=n\pi^2$ with $n=1,2,3,...$, and $[2n\pi]^2$ and $[2(n+1)\pi]^2$ are the limits of the domain of $n$-revolution elliptical orbits. Also in this case, there exists a lower limit for the value of $z$ for transfers with $\theta<\pi$. The work also provides some insight on how to avoid and correct certain scenarios where numerical precision gets degraded.

Battin and Vaughan \cite{BATTIN1984} improve the original Gauss formulation in their Elegant Algorithm, extending the convergence region and treating $\theta=\pi$ transfers. The new system of equations is also expressed in terms of the sector area-to-triangle area ratio $y$ and makes use of a continued fraction. The method, however, is singular for $2\pi$ transfers.

Gooding \cite{Gooding1990} follows the work of Lancaster \& Blanchard \cite{Blanchard1969}in terms of the universal variable $x$ ($x^2 = 1 - s/2a$) (FALTA DEFINIR s ?), implementing a Halley cubic iterator, initial guess value formulas to ensure fast convergence, a power-series expansion in the near-parabolic range to avoid round-off errors, and expanding its multi-revolution capabilities. In his paper, Gooding provides a Fortran code with three routines: one computes the time of flight (and several of its derivatives) as a function of $x$ and using the formulation by Lancaster \& Blanchard. The second routine finds $x$, while implementing all the initial guesses for all orbits and multi-revolution cases. The main code computes the geometric parameters of the transfer, calls the second routine to compute $x$ for each multi-revolution case, and computes the velocity vectors.

Izzo \cite{Izzo2015} also starts from the work by Lancaster \& Blanchard \cite{Blanchard1969}, implements a new initial guess based on empirical results and a House-Holder iteration scheme that converges in about two or three iterations. The algorithm has been implemented in a python code available as part of the ESA PyKEP toolbox [REF]. For values of $x$ close to one (i.e., near-parabolic range), Izzo makes use of the Lagrange and Battin equations in terms of the universal variable $z$. As part of the same toolbox, Izzo also implements another version of the algorithm using only the Lagrange equation, a set of transformed variables and a Regula-Falsi iteration scheme.

\subsection{Solution algorithm}
Several methods can be used to numerically solve the equations resulting from the several treatments described in Sect. \ref{sec:tof}. Each method has its benefits and shortcomings. For example, the Newton-Raphson method needs a single starter and has a high convergence rate but works well only with single-valued monotonic functions. Singularities or steep changes in the derivative slow down the convergence, make the method converge to the wrong solution, or even diverge. The bisection method, requires two starters bracketing the solution and it is usually slower than Newton-Raphson, but it will always converge to the the correct solution when properly bracketed. The iteration methods implemented for the selected Lambert algorithms are presented below.

The Newton-Raphson scheme is a Householder root-finding method of 1$^{st}$ order. The method requires the derivative of the transfer time with respect to the iteration variable, which we denote by $w$. At the $i^{th}$ iteration, $w$ is updated through:
\begin{equation}
w_{i+1} = w_i - \frac{t - TOF}{\frac{dt}{dw}}. \label{eq:newton-raphson}
\end{equation}
$TOF$ is the target value for the time of flight, $t$ is current value of the transfer time and $\frac{dt}{dw}$ is its derivative with respect to $w$. The Newton-Raphson method is used by Lagrange and Bate.

The bisection method does not require any derivatives, since it evaluates the transfer time equation at the arithmetic mean point (e.g., $\bar{w}$) of a solution-bracketing interval (namely $[w_{low},w_{up}]$). Then, the resulting value for the time of flight $t$ is compared against the target $TOF$, and the interval bounds are updated as follows:
\begin{eqnarray}
w_{low} & = & \bar{w} \quad \quad t \le TOF \\ \nonumber
w_{up} & = & \bar{w} \quad \quad t > TOF \label{eq:bisection}
\end{eqnarray}
The bisection method is implemented in the Bate and Sim\'o algorithms.

The Gauss algorithm uses a successive substitution technique called the Gauss Method. The scheme recomputes all parameters with the current value of the free variable and updates it using the recomputed value of the parameters. Battin also develops a successive substitution solving method, exploiting the property that the largest positive real root of the cubic equation is always the correct solution, and using the resulting equation to directly compute the solution for the cubic.

Gooding uses the Halley method, which is a Householder root-finding method of $2^{nd}$  order. The first and second derivatives of the transfer time equation with respect to the free parameter are then required and used as follows:
\begin{equation}
w_{n+1} = w_n + \frac{T dt}{dt^2 + T \frac{d^2t}{2}} \label{eq:gooding_halley}
\end{equation}
where $T = TOF - t$, and $dt$, $d^2t$ are the first and second order derivatives of the time equation, respectively. 

\cite{Izzo2015} uses a $3^{rd}$-degree Householder method to iterate on the transfer time equation:
\begin{equation}
w_{n+1} = w_n-T\frac{dt^2 - T d^2t/2}{dt (dt^2 - T d^2t) + d^3t d^2t/6} \label{eq:izzo_householder}
\end{equation}
where $T = TOF - t$, $dt$ and higher-order derivatives are computed using Lancaster\& Blanchard's version of the transfer time equation.

In the alternative version of the algorithm at the PyKEP toolbox, Izzo uses a Regula-Falsi scheme as iteration method. The Regula-Falsi method requires a solution-bracketing interval as a starter and updates the iteration variable as follows:
\begin{equation}
w_{n+1} = \frac{xw1 T_2 - T_1 w_2}{T_2 - T_1} \label{eq:izzo_regulafalsi}
\end{equation}
where $w_1$, $w_2$ are interval limits, while $T_1$ and $T_2$ are the values for the difference $t_n-TOF$ at each interval limit. Once the new value of the iterator ($w_{n+1}$) is computed, the corresponding value of the time of flight ($t_{n+1}$) is obtained via the transfer time equation. Then, the interval is updated as follows: $w_1 = w_2$, $w_2 = w_{n+1}$, $T_1 = T_2$ and $T_2 = T_{n+1}$. The process is repeated until convergence.

\subsection{Choice of the initial guess}
The initial guess, or starter, is introduced as input to the solution algorithm. Note that each solution algorithm may require a specific number of initial guesses (e.g., Householder-type methods require a single starter, bisection requires two starters, etc.). The starters proposed by the authors of the selected algorithms are the following:

Lagrange starts from $x = 0$, corresponding to a minimum-energy orbit.

Gauss starts from $y = 1$, which corresponds to a 1:1 area-to-sector ratio.

Bate starts from $z = 0$ for the Newton-Raphson scheme, corresponding to a parabolic orbit. For the bisection method, Vallado \cite{Vallado1997} suggests an upper bound of $z_{up} = 4\pi^2$, corresponding to the $t=2\pi$ asymptote, and a lower limit of $z_{low} = -4\pi$, which Vallado claims is valid for most orbits except for highly-eccentric ones. In such cases, the lower limit should be extended.

Sim\'o proves that his time equation has a single solution in the interval $z \in [z_f,\pi^2]$, where $\pi^2$ is the $t = 2\pi$ asymptote and $z_f$ is the lower limit for the $z$ parameter. The value of $z_f$ corresponds to a specific value for $\theta < \pi$ transfers and $z_f = -\infty$ for $\theta > \pi$ transfers. Sim\'o recommends to use $z = \left(\frac{\theta}{2} \right)^2$ as a starter for the elliptical orbits domain ($z \in [0,\pi^2]$) and $z = 0$ as a starter for hyperbolic orbits ($z \in [z_f,0]$).

Battin provides a specific initial guess $x$ for his successive substitution method if the normalized time of flight corresponds to an elliptical orbit, and $x = 0$ otherwise.

Gooding derives a set of semi-empirical initial starters for both single-revolution and multi-revolution transfers. For the single-revolution case, he uses a bilinear approximation method to analytically obtain an approximate solution for the $x<0$ region of Lancaster's transfer time equation, and a weighted combination of approximate solutions for the $x>0$ region. The multi-revolution starters are far more complex and, for reasons of space, the reader is referred to the corresponding publications.

Izzo starts performs a piece-wise linear approximation of the transfer time equation. Inverting the previous relation provides an approximate solution, which is then transformed back into the $x-T$ plane and used as a starter. For the multi-revolution cases, Izzo performs the same approximation with the lower and upper asymptotes of the time equation for the corresponding $n$-revolution scenario (?).
For his Regula-Falsi implementation, Izzo uses an empirical set of constant starters for the single-revolution case or multi-revolution case, both computed in the transformed $\xi-\tau$ plane. (DTS: arreglar)

\subsection{Computation of the velocity vectors}
Once the free parameter has been approximated, the velocity vectors at the boundaries of the transfer must be determined. The methods employed by the several authors are highlighted below.

\subsubsection{Orbital elements}
This method consists in expressing the semi-latus rectum $p$, the semi-major axis $a$, the eccentricity $e$ and the true anomaly $\nu$ of the selected point as functions of the free parameter. Then, the velocity ${\bf v}_{pwq}$ at either ${\bf r}_1$ and ${\bf r}_2$ in the perifocal reference frame is determined as:
\begin{equation}
{\bf v}_{pqw} = \sqrt{\frac{\mu}{p}} \left[-\sin \nu, e + \cos \nu, 0 \right].
\end{equation}
Then, with the inclination $i$, the argument of the periapsis $\omega$ and the longitude of the ascending node $\Omega$, the pqw vectors can be rotated into the xyz frame:
\begin{equation}
{\bf v}_{ijk} = R_z(\Omega) R_x(i) R_z(\omega) {\bf v}_{pqw},
\end{equation}
where $R_k$'s are the elementary rotation matrices following the right-hand sign convention. The orbital elements method is used by Sim\'o and Battin.

\subsubsection{Gauss $f$ and $g$ functions}
The values for the Gauss functions $f$, $g$ and $\dot{g}$ are expressed as functions of the free parameter. Then the velocity vectors ${\bf v}_1$ and ${\bf v}_2$ at the beginning and at the end of the transfer, respectively, are found through:
\begin{eqnarray}
{\bf v}_1 & = & ({\bf r}_2 - f {\bf r}_1) / g \\
{\bf v}_2 & = & (\dot{g} {\bf r}_2 - {\bf r}_1) / g 
\end{eqnarray}
Note that there exists a singularity in the computation of the velocity when $g = 0$, which corresponds to $\theta = \pi$ transfers. This method is used by Gauss and Bate.

\subsubsection{Radial and transversal components}
The radial $v_r$ and transversal $v_t$ components of the velocity are expressed as functions of the free parameter. Then, either velocity vector  ${\bf v}_i$ (i.e., ${\bf v}_1$ or ${\bf v}_2$) is calculated as:
\begin{equation}
{\bf v}_i = v_{r_i} \frac{{\bf r}_i}{r_i} + v_{t_i} \frac{{\bf h} \times {\bf r}_i}{|{\bf h} \times {\bf r}_i|}.
\end{equation}
Here, ${\bf r}_i/ r_i$ is the unit vector in the radial direction at $P_i$, and ${\bf h} \times {\bf r}_i /|{\bf h}\times {\bf r}_i|$ is the unit vector in the transversal direction. Lagrange, Gooding and Izzo use this method.

\section{Performance tests}
\label{sec:performance}
The methods discussed in Section~\ref{sec:methods} have been implemented following the description of their authors. Only single-revolution cases have been considered for this contribution. The implementation has been done in Matlab. This programming framework allows for a fast prototyping and includes powerful code analysis and debugging tools, which favours the algorithm-oriented analysis that we aim at in this work. Although a preliminary performance-oriented analysis has also been carried in order to obtain a relative comparison between methods, the correct practice in this case should be to programme the algorithms in a high-performance compiled code (e.g., FORTRAN or C). The latter approach will be addressed in a future development of this work.

\subsection{Application region}
\label{sec:region}
Several simulations have been carried out in order to test the algorithms in a wide range of scenarios. The goal is to stress each method and detect its range of applicability, slow-convergence zones, singularities, etc. The scenario sets a transfer between two dradial distances, $r_1$ and $r_2$, such that $r_2 / r_1 = 2$. Then the transfer angle $\theta$ is varied between $0$ and $2 \pi$ at constant steps, and for each value of the transfer angle, the transfer time is also varied between $0$ and $2 \pi$. The normalization has been applied which consists in that $\mu$ and $r_1$ are unitary and the circular orbit with radius $r_1$ has unit mean motion. This setup ensures that the algorithms are tried also in the tough cases, including $\theta = 0,\pi,2 \pi$, parabolic orbits, highly-eccentric hyperbolic orbits. The results of the simulations are discussed below.

\begin{figure}[ht]
\centering
\includegraphics[scale=0.9]{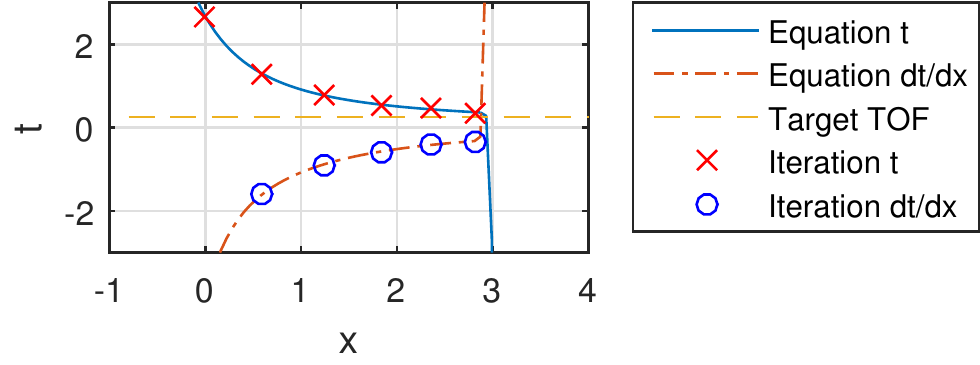}
\caption{Lagrange Newton-Raphson solver in the low-$TOF$ scenario.}
\label{fig:debug_lagrange_nr_01}
\end{figure}

\begin{figure}[ht]
\centering
\includegraphics[scale=0.75]{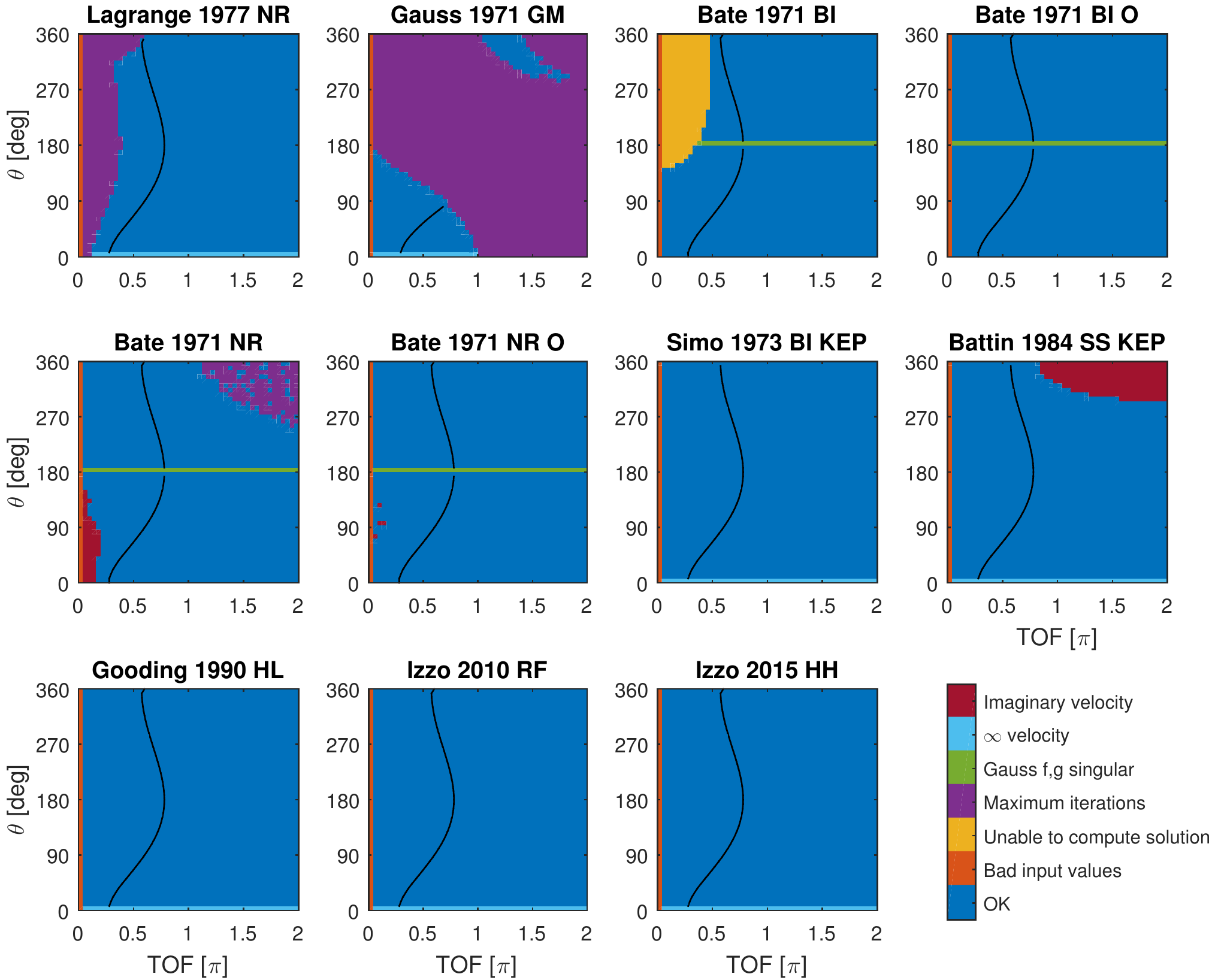}
\caption{Region of application and degeneracies of Lambert's algorithms.}
\label{fig:region_flag}
\end{figure}

Battin's implementation of Lagrange's equation fails for highly-eccentric hyperbolic orbits with $TOF < \frac{\pi}{4}$ (Fig.~\ref{fig:region_flag}), corresponding to values of the universal variable $x \gg 1$. In this scenario, the Gauss hypergeometric function $_2F_1$ does not converge because its convergence radius ($R=1$) is exceeded. As a result, the equation suffers a steep change in the derivative with which the Newton-Raphson method can no longer cope (Fig.~\ref{fig:debug_lagrange_nr_01}). The method performs 10-20 iterations in the elliptical orbits domain (best case near $\theta \approx \pi$), and executes 40-60 iterations when approaching the non-convergence zone of $_2F_1$.

The Gauss Method algorithm expanded by Bate only works for scenarios where $y \approx 1$ (equal area-to-sector ratio), which includes both elliptical orbits with a small transfer angle and hyperbolic orbits with a moderate transfer angle. In general, the convergence region is given by the approximate relation $\theta < \pi - TOF$, as shown in Fig.~\ref{fig:region_flag}. The method, however, fails completely in other cases. It is worth mentioning that for $\theta > \pi$, the value of the sector-to-triangle area ratio $y$ becomes negative, which corresponds to a physically impossible situation. Also, there is a small region for high values of $\theta$ where the algorithm converges but the results are obviously incorrect. The precision of the results in the convergence region is also relatively low compared with other algorithms (more than $10^{-2}$ relative error on the magnitude of the velocity vector). Due to the adoption of the $f$ and $g$ functions for the computation of the velocity vectors, the method also intrinsically fails for the $\theta = \pi$ transfer. The Gauss Method by Bate makes approximately 30 iterations in the elliptical orbits zone and 10-20 in the hyperbolic orbits domain, but the iteration number quickly increases to more than 100 when moving away from the convergence region.

\begin{figure}[ht]
\centering
\includegraphics[scale=0.9]{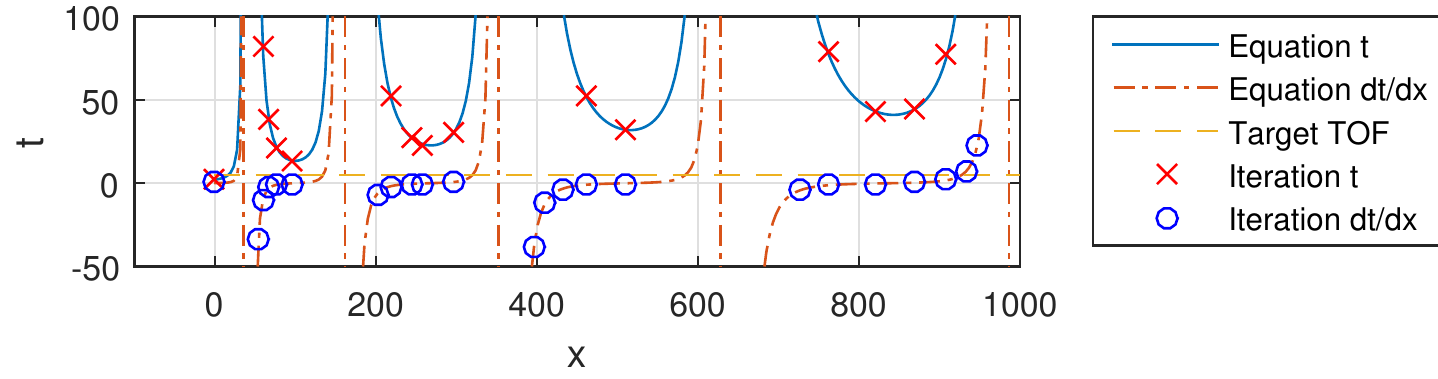}
\caption{Bate's Newton-Raphson solver for the scenario with high values of $TOF$ and $\theta$.}
\label{fig:debug_bate_nr_01}
\end{figure}

\begin{figure}[ht]
\centering
\includegraphics[scale=0.9]{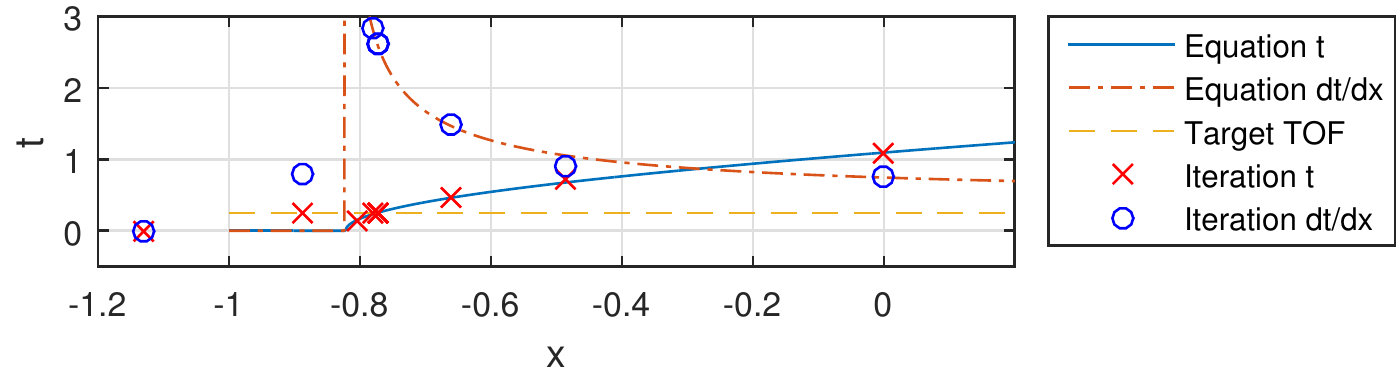}
\caption{Bate's Newton-Raphson solver for the scenario with very low $TOF$ and $\theta < \pi$.}
\label{fig:debug_bate_nr_02}
\end{figure}

Bate's method using universal variables and Newton-Raphson with a constant starter $z=0$ (parabolic orbit) has two conflictive regions. One corresponds to large values of $TOF$ and high $\theta \nu$ in the elliptical orbits range (Fig.~\ref{fig:region_flag}), where the solution is a point close to the $z=4 \pi^2$ asymptote and the derivative is steep (?). In this situation, Newton-Raphson has a high chance to "jump" over the asymptote and fall into the multi-revolution region, where the error is propagated further (Fig.~\ref{fig:debug_bate_nr_01}). The second conflictive zone is associated to $\theta < \pi$ (Type I, short-way transfers) and $TOF < \frac{\pi}{8}$ in the hyperbolic orbits region, where the solution falls close to the minimum negative value for $z$. In this case, Newton-Raphson has a high chance of falling into the imaginary region of Bate's equation (Fig.~\ref{fig:debug_bate_nr_02}), thus yielding imaginary numbers. These two failure regions can be avoided by improving the initial guess, for example with a preliminary search to bracket the solution within bounds and then selecting the geometric mean between bounds as a starter. In this case, the convergence is greatly improved, except for a few cases with extremely low or extremely high values of $TOF$. The number of iterations in any case ranges from 5-10 starting from the parabolic orbits region.  [DTS arreglar esta frase]

Bate's method with bisection and the starter proposed by Vallado works for all orbits except those with $\theta > \frac{3}{4} \pi$ and $TOF < \frac{\pi}{2}$ in the hyperbolic orbits region, where the initial bounds no longer bracket the solution. In such cases, Vallado suggests to decrease the value of the lower starter until the solution is bracketed again. Of course, this new search implies a slight increase in computing time, but convergence is eventually achieved in all regions. The iteration counter is fairly constant and close to 40 iterations. Note that both Bate's methods fail to compute the velocity vectors for $\theta = \pi$ due to the singularities in the $f$ and $g$ functions.

Sim\'o's algorithm presents the same conflictive regions as Bate's, as expected. The improvement comes again from a more careful choice of the starter, which in this case has been implemented also with a preliminary search for a solution-bracketing bounds. With the bisection method, the convergence is guaranteed in all regions, as long as the solution is properly bracketed first. Sim\'o computes the velocity vectors with the method of the orbital elements, which gives good results for all orbits (including $\theta = \pi$ transfers). The iteration counter stays fairly constant around 40 iterations in the $0<TOF<\pi$ region and close to 35 iterations in the $\pi<TOF<2\pi$ region.

Battin's Elegant Algorithm implemented using Battin's successive substitution method with direct computation of the cubic root fails for large values of $TOF$ and $\theta$, as illustrated in Fig.~\ref{fig:region_flag}. In this scenario, one of the parameters in Battin's equations produces a negative square root and therefore the method gives imaginary solutions. The algorithm converges in all other regions, although the precision is relatively low ($10^{-5}-10^{-3}$ errors relative to the other algorithms). The algorithm converges in about 4-5 iterations for the region $\theta < \pi$ and $TOF < \pi$ (the same region as Gauss), 8-9 iterations for hyperbolic orbits with $\theta > \frac{3}{2} \pi$, and 6-7 iterations elsewhere.

Gooding's algorithm has no evident singularities in the entire single-revolution region, except when $TOF = 0$ and $\theta = 0$ which give infinite velocity. The algorithm performes a fixed number of iterations, as Gooding deminstrates that three iterations are sufficient to provide 13 digits precision for all cases in a single-revolution scenario. However, the error slowly increases with the multi-revolution number, which would require the convergence to be re-assessed.

Both versions of Izzo's algorithm converge in all cases except for $TOF = 0$ and $\theta = 0$. The Regula-Falsi method converges in 5-6 iterations, while the Householder version converges in just 2-3 iterations.

A combined view of the regions of application of the several algorithms is shown in Fig.~\ref{fig:region_flag}. The -O (DTS arreglar) versions of both Bate algorithms implement an improved initial guess estimation.

\subsection{Computational performance}
A computational performance test has been carried out in order to estimate the execution time of each algorithm. The algorithms have been tested on a Monte-Carlo simulation of $10^5$ random distributed samples within the single-revolution region defined in Section~\ref{sec:region}, and only the points for which the algorithm actually converges correctly have been used to compute the resulting average execution time. This procedure aims at simulating a realistic application scenario.

While implementing the algorithms, extremely large improvements in execution time were observed whenever the code was slightly optimised. Certainly, as the algorithms take of the order of microseconds to execute, the normally small overhead caused by callbacks to external functions (e.g., norm, dot, cross, Stumpff functions, series, etc.) and inefficient computations (e.g., repeating trigonometric calculations) end up having a great impact on the overall execution time. In the end, all the codes were optimized by using only intrinsic functions, pre-computing repeated calculations, minimising the number of calls to external functions and subroutines whenever possible, and keeping a strictly uniform coding style among the implementations. Note that any different implementation will most probably produce slightly different results. 
%As Klumpp mentioned in his report \cite{Klumpp1991a}: ``\textit{...every [code] implementation is unique}''.

\begin{figure}[ht]
\centering
\includegraphics[scale=0.75]{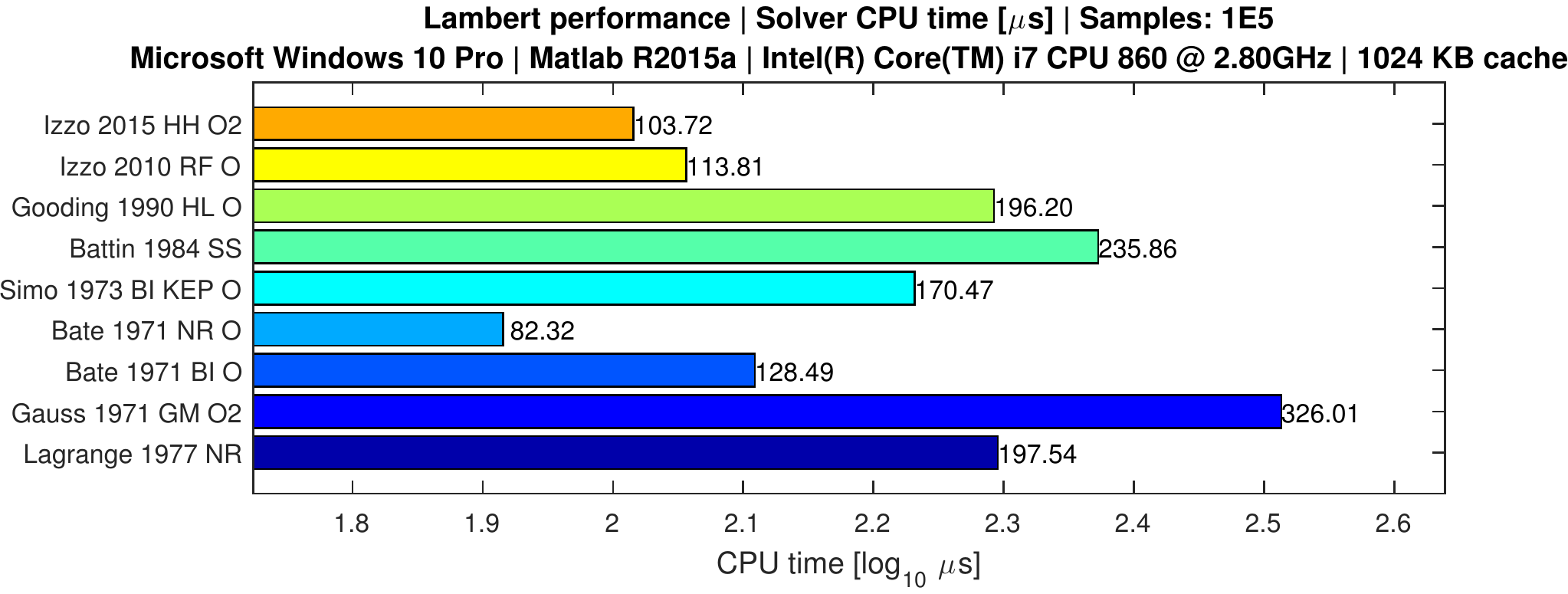}
\caption{CPU execution time of Lambert algorithms}
\label{fig:cpu_time}
\end{figure}

The results of the performance comparison are presented in Fig.~\ref{fig:cpu_time}. The reader should take these results with care: first of all, the computing speed of Matlab's interpreted code is largely slower than compiled code (around 8 times slower). The Matlab processor may also have automatically used unexpected threading and optimization capabilities on specific portions of the code that could eventually alter the results. Finally, executing the code under a UNIX system would have provided better results in terms of speed. Hence, these values can be taken as a first-order approximation in a comparison among the algorithms, or if the application is to be run in the Matlab platform.

\section{Conclusions}
\label{sec:conclusions}
The existing bibliography on Lambert's problem has been reviewed and organized. The major existing algorithms have been selected, analysed  and compared in terms of formulation of the transfer time equation, initial guess estimation, choice of the free parameter, iteration method, and velocity vector computation scheme. A series of tests have been carried out to determine the range of applicability and the computational performance of each algorithm. Bate's algorithm via universal variables with Newton-Raphson appears to be the fastest solver, but Izzo's Householder algorithm exhibits the best ratio between speed, robustness and accuracy.

The foreseen developments include implementing the algorithms in compiled code (FORTRAN), testing the implementation under a GPU computing framework (CUDA), incorporating additional algorithms, and expanding the analysis to include multi-revolution solutions.
\section*{Acknowledgements} 
David de la Torre acknowledges financial support by DLR.
The authors are grateful to Prof. Gerard G\'omez for suggesting this research and for helpful discussions throughout.

\bibliographystyle{ISSFD_v01}
\bibliography{Lambert_Bibliography}

\end{document}